\documentclass[letterpaper,11pt]{article}

\title{A Finiteness Theorem for Markov Bases of Hierarchical Models}

\author{Serkan Ho{\c s}ten  \\
{\small Department of Mathematics,
San Francisco State University, San Francisco} \\
\and Seth Sullivant
 \\
{\small Society of Fellows, Harvard University} }

\date{}

\usepackage{latexsym,array,delarray,amsthm,amssymb,epsfig}
\usepackage[sumlimits]{amsmath}

\theoremstyle{plain}
\newtheorem{thm}{Theorem}[section]
\newtheorem{lemma}[thm]{Lemma}

\newtheorem{cor}[thm]{Corollary}
\newtheorem{conj}[thm]{Conjecture}

\theoremstyle{definition}
\newtheorem{defn}[thm]{Definition}
\newtheorem{ex}[thm]{Example}

\theoremstyle{remark}
\newtheorem*{rmk}{Remark}

\newcommand{\zz}{\mathbb{Z}}
\newcommand{\nn}{\mathbb{N}}

\newcommand{\bft}{\mathbf{t}}
\newcommand{\bfu}{\mathbf{u}}
\newcommand{\bfv}{\mathbf{v}}

\begin{document}
\maketitle

\begin{abstract}
We show that the complexity of the Markov bases of
multidimensional tables stabilizes eventually if a single table dimension is
allowed to vary.  In particular, if this table dimension
is greater than a computable bound, the Markov bases consist of
elements from Markov bases of smaller tables.
We give an explicit formula for this bound in terms of Graver bases.
We also compute these Markov and Graver complexities for
all $K \times 2 \times 2 \times 2 $ tables.
\end{abstract}


\section{Introduction}

Let $d_1, \ldots, d_n$ be positive integers where $d_i \geq 2$. A
multidimensional contingency table is an $d_1 \times \ldots \times
d_n$ array of nonnegative integers. Such a table represents the
results of a census of individuals for which $n$ discrete random
variables $X_1, \ldots, X_n$ are observed (where we assume the random variable
$X_i$ takes values in $[d_i] := \{1, \ldots, d_i\}$) .
Inferences about the collected data are made based on a 
statistical model or collection of models. This paper is concerned with
the family of {\em hierarchical log-linear models} for which one
assumes a set of interaction factors between the
random variables \cite{F}.  Performing the exact test of conditional inference requires knowledge of the Markov basis of a given hierarchical model, which we describe below.

When we assume that the sampling
distribution of a table of observations is Poisson or multinomial,
the sufficient statistics of any hierarchical model are given by
certain marginal totals.  The particular marginal totals that are
sufficient statistics depend on the hierarchical model. For
instance, for a $d_1 \times d_2 \times d_3$ contingency
table, the \emph{no three-way interaction model} has sufficient
statistics that are the three 2-way  margins of the table:
$$ u_{+jk} = \sum_{i =1}^{d_1} u_{ijk} , \,\,\,\,\,\,
 u_{i+k} = \sum_{j =1}^{d_2} u_{ijk} , \,\,\,\,\,\,\,
 u_{ij+} = \sum_{k =1}^{d_3} u_{ijk}, $$
where $u_{ijk}$ are the entries of the table.

In general,  a hierarchical model (and hence the marginal totals)
is described by the list of the maximal faces $F_1, \ldots, F_r$
of a  simplicial complex $\Delta$ on $n$ vertices. Computing
marginal totals corresponds to a linear map from the space of tables 
to the space of marginals:
$$ \pi_\Delta \,\, : \,\, \nn^{D} \longrightarrow \bigoplus_{k=1}^r \nn^{D_k}$$
where $D = \prod_{i=1}^n d_i$ and $D_k = \prod_{j \in F_k} d_j$.
This map is defined by
$$(u_{i_1, \ldots, i_n}: \,\, i_j \in [d_j]) \longrightarrow
  \bigoplus_{k = 1}^r \sum_{(i_j: \, j \notin F_k)} u_{i_1, \ldots, i_n}.$$
Two tables $\bft$ and $\bfu$ are said to be in the same \emph{fiber} of
$\pi_\Delta$ if $\pi_\Delta(\bft) = \pi_\Delta(\bfu)$. In other words,
two tables are in the same fiber if they have the same margins with respect to $\Delta$. We
say that $\bft$ and $\bfu$ are connected by the sequence of \emph{moves} 
$\bfv_1, \ldots, \bfv_s$ if each move $\bfv_i$ is in  
$\ker_\zz(\pi_\Delta)$ 
(that is, $\bfv_i$ has zero margins), $\bft + \sum_{i=1}^p \bfv_i$ is a 
table with
nonnegative entries for each $1 \leq p \leq s$, and  $\bfu = \bft +
\sum_{i=1}^s 
\bfv_i$. 

By a theorem of Diaconis and Sturmfels \cite[Theorem 3.1]{DS}, for each
hierarchical model given by $\Delta$ and $d = (d_1, \ldots, d_n)$,
there exists a finite set of moves called a \emph{Markov basis}
such that any
two tables that are in the same fiber of $\pi_\Delta$ are connected
by the moves in the Markov basis.  Computing Markov bases
via \emph{Gr\"obner bases} for the use in MCMC methods was initiated
in \cite{DS}, and since this first work computing Markov bases
efficiently and describing Markov bases succinctly have been the major focuses of research.
Recently substantial progress has been made. Simple Markov bases
(consisting of moves with four nonzero entries)  for {\em decomposable
models} have been determined \cite{DF}, and similar
Markov bases are known for {\em reducible models} \cite{DoSu, HS}.
The case of \emph{binary graph models} (where $d_1 = \cdots = d_n = 2$
and $\Delta$ is a graph on $n$ vertices) is worked out up to $n=5$ \cite{DeSu}.

The contribution of this paper is the most general form of a
result first obtained in \cite{AT} for the no three-way
interaction model for $K \times 3 \times 3$ tables.  Our main
theorem and its proof rely on ideas from \cite{SS} which treats
the case of $K \times d_2 \times \cdots \times d_n$ tables where
$\{2,3,\ldots, n \}$ is a maximal face of  $\Delta$-- the
so-called {\em logit models}.

\begin{thm} \label{thm:main}
Let $\Delta$ and $d = (d_1, \ldots, d_n)$ define a hierarchical model.
Then there exists
a constant $ m := m(\Delta; d_2, \ldots, d_n)$
such that for all $d_1 \geq m$ the  universal Markov basis $\mathcal{M}_{\Delta, d}$
consists of tables of the format
$r \times d_2 \times \cdots \times d_n$ where $r \leq m$.
\end{thm}

In other words, if we fix a hierarchical model together with the
$n-1$ dimensions $d_2, \ldots, d_n$ while varying the single
dimension $d_1$, then for large enough $d_1$ the universal Markov basis
$\mathcal{M}_{\Delta, d}$ will be obtained from Markov bases of
small fixed-size tables; see Definition \ref{defn:markov}.  
In particular, as a function of 
$d_1$, the complexity of computing and storing a Markov basis is bounded.

We give the details of the proof of
Theorem \ref{thm:main} in Sections 2 and 3, where we give an explicit
computable upper bound for $m(\Delta; d_2, \ldots, d_n)$.  We also
present, in Section 3, a lower bound for $m(\Delta; d_2, \ldots,
d_n)$ which only applies in some cases.  In the fourth section, we consider strengthenings and generalizations of the main result when more than one level is allowed to vary.
In the final section we explicitly compute
the complexity bound $m(\Delta;2,2,2)$ for all tables of the
form $K \times 2 \times 2 \times 2$.


\section{From Models to Matrices}

In this section  we describe how to obtain a matrix $A_\Delta$
corresponding to the linear transformation $\pi_\Delta$ where
the simplicial complex $\Delta$ describes a hierarchical
model. We describe a decomposition for $A_\Delta$ that
is fundamental to the proof of Theorem \ref{thm:main}.

Given the vector  $d = (d_1, \ldots, d_n)$, and a
subset $F = \{j_1 < \ldots < j_s\}  \subseteq [n]$
we let $d_F = (d_{j_1},d_{j_2}, \ldots, d_{j_s})$.
The columns of $A_\Delta$ are in bijection with the $D$ entries
of a $d_1 \times \cdots \times d_n$ table, and we label each
such column with the vector indexing the table entry
$(i_1, \ldots, i_n) \in [d_1] \times \cdots \times [d_n]$.
Moreover, we order these columns lexicographically:
$$(1,1, \ldots,1, 1) \prec (1,1, \ldots,1,2) \prec \cdots \prec
(1,1, \ldots, 1, d_n) \prec $$
$$(1,1,\ldots, 2,1) \prec  (1,1,\ldots, 2,2) \prec \cdots  \prec
(1,1, \ldots, 2, d_n) \prec \cdots $$
$$(d_1, d_2, \ldots, d_{n-1}, 1) \prec  (d_1, d_2, \ldots, d_{n-1}, 2)
\prec  \cdots \prec (d_1, d_2, \ldots, d_{n-1}, d_n).$$
Each row is labeled by a pair $(F,e)$ where
$F = \{j_1, j_2, \ldots, j_s\}$  is a facet of
$\Delta$ and $e = (e_{j_1}, e_{j_2}, \ldots, e_{j_s}) \in
[d_{j_1}] \times [d_{j_2}] \times \cdots \times [d_{j_s}]$ indexing the
marginal corresponding to $F$. We first list the rows $(F,e)$ where $1 \in F$.
We impose a linear order on the facets where $1 \in F$ and set $(F,e) 
\prec (G,f)$
if $e_1 < f_1$, or if $e_1 = f_1$ and $F \prec G$, and in the case
when $e_1=f_1$ and $F=G$ we use an arbitrary but fixed order of the
indices. The rest of the rows will be listed again by some arbitrary
but fixed order which will not play a role for the rest of the article.
The entry of $A_\Delta$ in the column indexed by $(i_1, \ldots, i_n)$
and the row $(F =\{j_1, \ldots, j_s\} , (e_{j_1}, \ldots, e_{j_s}))$ will
be equal to one if $i_{j_1} = e_{j_1}$, $i_{j_2} = e_{j_2}, \ldots$, and
$i_{j_s} = e_{j_s}$; and it will be zero otherwise.
\begin{ex}\label{ex:model}
Let $\Delta = \{\{1,2\},\{1,4\},\{2,3\}, \{3,4\}\}$ and $d_1=d_2=d_3=d_4=2$.
This is the binary $4$-cycle model.
\renewcommand\arraystretch{0.8}
\fontsize{10}{12}
$$
A_\Delta \quad = \quad \left[ \begin{array}{cccccccc|cccccccc}
 1 &  1 &  1 &  1 &  0 &  0 &  0 &  0 &  0 &  0 &  0 &  0 &  0 &  0 &  0 &  0 \\
 0 &  0 &  0 &  0 &  1 &  1 &  1 &  1 &  0 &  0 &  0 &  0 &  0 &  0 &  0 &  0 \\
 1 &  0 &  1 &  0 &  1 &  0 &  1 &  0 &  0 &  0 &  0 &  0 &  0 &  0 &  0 &  0 \\
 0 &  1 &  0 &  1 &  0 &  1 &  0 &  1 &  0 &  0 &  0 &  0 &  0 &  0 &  0 &  0 \\
\hline
 0 &  0 &  0 &  0 &  0 &  0 &  0 &  0 &  1 &  1 &  1 &  1 &  0 &  0 &  0 &  0 \\
 0 &  0 &  0 &  0 &  0 &  0 &  0 &  0 &  0 &  0 &  0 &  0 &  1 &  1 &  1 &  1 \\
 0 &  0 &  0 &  0 &  0 &  0 &  0 &  0 &  1 &  0 &  1 &  0 &  1 &  0 &  1 &  0 \\
 0 &  0 &  0 &  0 &  0 &  0 &  0 &  0 &  0 &  1 &  0 &  1 &  0 &  1 &  0 &  1 \\
\hline
 1 &  1 &  0 &  0 &  0 &  0 &  0 &  0 &  1 &  1 &  0 &  0 &  0 &  0 &  0 &  0 \\
 0 &  0 &  1 &  1 &  0 &  0 &  0 &  0 &  0 &  0 &  1 &  1 &  0 &  0 &  0 &  0 \\
 0 &  0 &  0 &  0 &  1 &  1 &  0 &  0 &  0 &  0 &  0 &  0 &  1 &  1 &  0 &  0 \\
 0 &  0 &  0 &  0 &  0 &  0 &  1 &  1 &  0 &  0 &  0 &  0 &  0 &  0 &  1 &  1 \\
 1 &  0 &  0 &  0 &  1 &  0 &  0 &  0 &  1 &  0 &  0 &  0 &  1 &  0 &  0 &  0 \\
 0 &  1 &  0 &  0 &  0 &  1 &  0 &  0 &  0 &  1 &  0 &  0 &  0 &  1 &  0 &  0 \\
 0 &  0 &  1 &  0 &  0 &  0 &  1 &  0 &  0 &  0 &  1 &  0 &  0 &  0 &  1 &  0 \\
 0 &  0 &  0 &  1 &  0 &  0 &  0 &  1 &  0 &  0 &  0 &  1 &  0 &  0 &  0 &  1
\end{array} \right]  $$
\normalsize
Here the 16 columns are indexed as
$(1,1,1,1) \prec (1,1,1,2) \prec (1,1,2,1) \prec (1,1,2,2) \prec \cdots
\prec (2,2,2,1) \prec (2,2,2,2)$. The first four rows are indexed by
$(F_1, (1,1)), (F_1, (1,2)), (F_2, (1,1)),$ and $(F_2, (1,2))$ where
$F_1 = \{1,2\}$ and $F_2 = \{1,4\}$. The second block of four rows are
indexed by $(F_1, (2,1)), (F_1, (2,2)), (F_2, (2,1)),$ and $(F_2, (2,2))$.
For the rest of the rows we have chosen the order $(F, (i,j)) \prec
(G, (s,t))$ if $F = \{2,3\}$ and $G= \{3,4\}$, or
if $F=G$ and $(i,j) \prec (s,t)$ lexicographically.
\end{ex}

We observe that when the rows and columns are ordered as described,  the matrix $A_\Delta$
exhibits a block structure. In the above example, the upper-left
and lower-right blocks of the first eight rows are identical augmented
with the two blocks of zeros. And the last eight rows are split into
two identical matrices. We summarize this observation in the
following lemma where we assume the ordering of the columns and
rows of $A_\Delta$ that we introduced above.
\begin{lemma}\label{lem:bloc}
Let $\Delta = \{F_1, \ldots, F_r\}$ and $d = (d_1, \ldots, d_n)$ define
a hierarchical model. Then
$$A_\Delta \quad = \quad
\left[ \begin{array}{ccccc}
A & 0 & 0 & \cdots & 0 \\
0 & A & 0 & \cdots & 0 \\
\vdots & \vdots & \ddots & \vdots & \vdots \\
0 & 0 & 0 & \cdots & A \\
B & B & B & \cdots & B
\end{array} \right]
$$
where $A$ is a $\sum_{k=1}^s (D_k/d_1)  \, \times \, (D/d_1)$ matrix
with $F_1, \ldots, F_s$ being the facets containing the vertex
$1$, and where $B$ is a $\sum_{k=s+1}^r D_k \, \times \, (D/d_1)$ matrix.
Hence there are $d_1$ copies of $A$ and $B$.
\end{lemma}

\begin{rmk}
The matrices $A$ and $B$ are also matrices that come from hierarchical models.
Note that the matrix $A$ is the matrix $A_{\Gamma}$ for the simplicial complex 
$$\Gamma = {\rm link}(\Delta) := {\rm link}(\Delta, \{1\}) = \{ F 
\setminus \{1\} \, | \, F \in \Delta \mbox{ and } 1 \in F \}$$
and the vector $d' = (2,3,\ldots, n)$.  The matrix $B$ is the matrix $A_{\Delta \setminus \{ 1\} } $ for the simplicial complex
$$\Delta \setminus \{1 \} = \{ F \, | \, F \in \Delta, 1 \notin F, \mbox{ and } F \mbox{ is a facet} \}.$$
\end{rmk}


\section{Proof of the Finiteness Theorem}

In this section, we proceed with the proof of Theorem \ref{thm:main}.  To do this, we will prove a finiteness theorem for the Markov bases of arbitrary matrices which come in a block form akin to the one demonstrated in Lemma \ref{lem:bloc}.

\begin{defn} \label{defn:markov}
Let $A \in \nn^{d \times n}$ be an integer matrix with no zero columns.  A finite set $\mathcal{M} \subset \ker_\zz(A)$ of integer vectors in the kernel of
$A$ is called a \emph{Markov basis}  of $A$ if any two nonnegative integer vectors
in the same fiber of $A$ can be connected by a collection
of the elements in $\mathcal{M}$.  That is, for any $\bft, \bfu \in \nn^n$ with $A\bft = A\bfu$, there exists a sequence of moves $\{\bfv_i \}_{i =1}^s \subset \mathcal{M}$ such that
$$ \bft  + \sum_{i = 1}^p \bfv_i \geq 0 \quad \mbox{for all } 1 \leq p \leq s \quad \mbox{ and }  \quad
\bft  + \sum_{i = 1}^s \bfv_i = \bfu.$$
A Markov basis $\mathcal{M}$ of $A$ is called \emph{minimal} if no subset of $\mathcal{M}$ is a Markov basis of $A$.  The \emph{universal Markov basis} $\mathcal{M}(A)$ of $A$ is the union of all minimal Markov bases of $A$.
\end{defn}

When $A = A_\Delta$ is the matrix associated to a hierarchical model, we use the shorthand $\mathcal{M}_{\Delta,d}$ to denote the universal Markov basis of $A_\Delta$.

\begin{defn} Let $\bfu, \bfv,$ and $\bfv'$ be nonzero vectors in $\ker_\zz(A)$.
We say $\bfu = \bfv + \bfv'$
is a \emph{conformal decomposition} of $\bfu$ if $u_i \geq 0$ implies
$0 \leq v_i, v_i' \leq u_i$, and $u_i \leq 0$ implies
$u_i \leq v_i, v_i' \leq 0$ for all $1 \leq i \leq n$. The
set $\mathcal{G}(A) \subset \ker_\zz(A)$ of integer vectors with no conformal decompositions
is called the \emph{Graver basis} of $A$.
\end{defn}

One can show that $\mathcal{G}(A)$ is a finite set \cite[Chapter 4]{GB+CP}
and any minimal Markov basis of $A$ is a subset of $\mathcal{G}(A)$
\cite[Chapter 5]{GB+CP}. Thus, $\mathcal{M}(A)$ is a finite set.

\begin{defn}
Let $A$ be a $d \times n$ matrix with columns $a_1, \ldots, a_n$ and $B$ be a $p \times n$ matrix with columns $b_1, \ldots, b_n$.  
The {\em $r$-th generalized Lawrence lifting} of $A$ with $B$ is the 
$(rd+p) \times rn$ matrix $\Lambda(A,B,r)$, whose columns are the vectors
$$\Lambda(A,B,r) = \{ a_i \otimes e_j \oplus b_i \, | \, 1 \leq i \leq n, 1 \leq j \leq r \}.$$
\end{defn}

In particular, the matrices $A_\Delta$ are of the form $\Lambda(A,B,r)$ with $A = A_{{\rm link}(\Delta)}$, $B = A_{\Delta \setminus \{1\}}$ and $r = d_1$.  When $B$ is the $n \times n$ identity matrix and $r=2$, the matrix
$\Lambda(A,B,r)$ is called the \emph{Lawrence lifting} of $A$
\cite[Chapter 7]{GB+CP}. For $B=I_n$ but general $r$, this matrix is
called the $r$th Lawrence lifting of $A$ \cite{SS}.

\begin{rmk} An integer vector in the kernel of $\Lambda(A,B,r)$
can be represented as an $r \times n$ matrix where each row is
in the kernel of $A$, and the sum of the rows is in the kernel of
$B$.  For instance, the following $2 \times 8$ matrix
is the representation of such a  vector in the kernel of
$A_\Delta$ in Example \ref{ex:model}:
$$\left[ \begin{array}{rrrrrrrr}
2 & -2 & -1 & 1 & -2 & 2 & 1 & -1 \\
-1 & 1 & 1 & -1 & 1 & -1 & -1 & 1
\end{array} \right].
$$
\end{rmk}

\begin{defn}
The \emph{type} of a vector in $\zz^{rn}$ represented as an
$r \times n$ matrix is the number of nonzero rows of this matrix. 
The \emph{Markov complexity} $m(A,B)$ of a $d \times n$ matrix $A$
and a $p \times n$ matrix $B$ is the largest type of any
vector in the universal Markov basis of $\Lambda(A,B,r)$ as $r$ varies.
Similarly, the \emph{Graver complexity} $g(A,B)$ of these two matrices
is defined as the largest type of any Graver basis element of 
$\Lambda(A, B,r)$ as $r$ varies. Analogously, we define $g(\Delta; d_2, 
\ldots, d_n)$ and $m(\Delta; d_2, \ldots, d_n)$, the Graver and Markov
complexities of the hierarchical models corresponding to $\Delta$ as
$d_1$ varies.
\end{defn}

We will show that the Graver complexity
$g(A,B)$ is finite. This implies that the Markov complexity
is also finite since $m(A,B) \leq g(A,B)$. In order to do this we
relate the Graver basis  $\mathcal{G}(B \cdot \mathcal{G}(A))$ to the collection of Graver
bases $\mathcal{G}(\Lambda(A,B, r))$. We emphasize the ``double'' Graver construction: we first compute
the Graver basis of $A$ and obtain the set $\mathcal{G}(A)$. We consider each element in $\mathcal{G}(A)$ as a column vector. Then
the vectors $B \cdot \mathcal{G}(A)$ are computed by multiplying each
element of $\mathcal{G}(A)$ with $B$. Thus, 
$B \cdot \mathcal{G}(A)$ is a $p \times |\mathcal{G}(A)|$ matrix. 
Finally we compute the Graver basis of $B \cdot \mathcal{G}(A)$.

\begin{thm} \label{thm:graver}
The Graver complexity $g(A,B)$ is the maximum 1-norm
of any element in the Graver basis  $\mathcal{G}(B \cdot \mathcal{G}(A))$.
\end{thm}

\noindent
In order to prove the above theorem we need the following lemma.

\begin{lemma} \label{lem:decomp}
Let $\bfu = [u^1; u^2; \ldots ; u^r]$ be in the Graver basis of
$\Lambda(A,B,r)$. Suppose that $u^i = v^1 + v^2$ is a conformal
decomposition where $v^1$ and $ v^2$ are in the kernel of $A$. Then
the element $[u^1; \ldots ; u^{i-1} ; v^1 ; v^2; u^{i+1} ; \ldots ; u^r]$
is in the Graver basis  of $\Lambda(A,B,r+1)$.
\end{lemma}
\begin{proof} Suppose not. Then
$[u^1; \ldots ; u^{i-1} ; v^1 ; v^2; u^{i+1} ; \ldots ; u^r]$ has
a conformal decomposition
$$[\bar{u}^1; \ldots ; \bar{u}^{i-1} ; \bar{v}^1 ; \bar{v}^2; \bar{u}^{i+1} ; \ldots ; \bar{u}^r] +
[\hat{u}^1; \ldots ; \hat{u}^{i-1} ; \hat{v}^1 ; \hat{v}^2; \hat{u}^{i+1} ; \ldots ; \hat{u}^r]$$ where both vectors are in the kernel of
$\Lambda(A, B,r+1)$. Now since $u^i = v^1 + v^2 = (\bar{v}^1 + \hat{v}^1) 
+
(\bar{v}^2 + \hat{v}^2)$ is a conformal decomposition of $u^i$, so
is $(\bar{v}^1 + \bar{v}^2) + (\hat{v}^1 + \hat{v}^2)$. We note that neither
the first nor the second sum is zero. But then the two nonzero
vectors
$[\bar{u}^1; \ldots ; \bar{u}^{i-1} ; \bar{v}^1 + \bar{v}^2; \bar{u}^{i+1} ; \ldots ; \bar{u}^r]$ and
$[\hat{u}^1; \ldots ; \hat{u}^{i-1} ; \hat{v}^1 + \hat{v}^2; \hat{u}^{i+1} ; \ldots ; \hat{u}^r]$ are in the kernel of $\Lambda(A,B,r)$, and their sum forms
a conformal decomposition of $\bfu$. This is a contradiction since $\bfu$ is in
the Graver basis of $\Lambda(A,B,r)$.
\end{proof}

\noindent \emph{Proof of Theorem \ref{thm:graver}}:  Lemma \ref{lem:decomp}
implies that in order to compute the Graver complexity $g(A,B)$ we only
need to consider elements $\bfu = [u^1; \ldots; u^r]$ where
$u^i$ is in the Graver basis $\mathcal{G}(A) = \{v_1,
\ldots, v_k\}$. Given any such $\bfu$, we construct a vector $\Gamma \in \zz^k$ where
the $i$th entry counts how many times $v_i$ appears in $\bfu$.
The $1$-norm of $\Gamma$ is the type of $\bfu$.  Hence we need
to show that $\Gamma$ is in the Graver basis of $B \cdot \mathcal{G}(A)$
if and only if $\bfu$ is in the Graver basis of $\Lambda(A,B,r)$.
If $\Gamma$ is not in the Graver basis of $B \cdot \mathcal{G}(A)$, then
it has a conformal decomposition $\Gamma_1 + \Gamma_2$ such
that $B \cdot \mathcal{G}(A) \cdot \Gamma_i = 0$ for $i=1,2$.  Reversing the operation that produced $\Gamma$ from $\bfu$, $\Gamma_1$ and $\Gamma_2$ yield vectors $\bfv_1$,  $\bfv_2 \in \ker_\zz(\Gamma(A,B,r)$ such that $\bfu = \bfv_1 + \bfv_2$ and this decomposition is conformal.
Thus $\bfu$ could not be in the Graver basis of $\Lambda(A,B,r)$.

Conversely, a conformal decomposition of $\bfu$ translates into 
a conformal decomposition of $\Gamma$ since none of  $u^1, \ldots, u^r \in \mathcal{G}(A)$
have conformal decompositions. \hfill $\Box$

\medskip

\noindent \emph{Proof of Theorem \ref{thm:main}}:  The hierarchical
model defined by $\Delta$ and $d$ gives
rise to $A_\Delta$ which is of the form $\Lambda(A_{{\rm link}( \Delta)},A_{\Delta \setminus \{1\}}, d_1)$ by Lemma \ref{lem:bloc}.  Theorem \ref{thm:graver} implies that
the Markov complexity $m(A_{{\rm link}( \Delta)},A_{\Delta \setminus \{1\}})$ is bounded by the finite
Graver complexity $g(A_{{\rm link}( \Delta)},A_{\Delta \setminus \{1\}})$. This means that for all 
$d_1 \geq m(A_{{\rm link}( \Delta)},A_{\Delta \setminus \{1\}})$
the universal Markov basis $\mathcal{M}_{\Delta, d}$ will consist
of tables of the format $r \times d_2 \times \cdots \times d_n$
where $r \leq m(A_{{\rm link}( \Delta)},A_{\Delta \setminus \{1\}})$. \hfill $\Box$ 

\medskip

In practice, the Graver complexity and the Markov complexity may vary a lot, 
as the following examples illustrate.  All of our examples were computed 
using {\tt 4ti2} \cite{He} and the results for Markov bases of reducible 
models using \cite{DoSu,HS}.

\begin{ex}
Let $\Delta = \{\{1,2\},\{1,3\},\{2,3\}\}$.  Then for $d_2 = d_3 = 3$, the Markov complexity is $m(\Delta; 3,3) = 5$ (the main result in \cite{AT}), while the Graver complexity is $g(\Delta;3,3) = 9$.
\end{ex}

\begin{figure}[h]\label{fig:triang}
\begin{center}\includegraphics{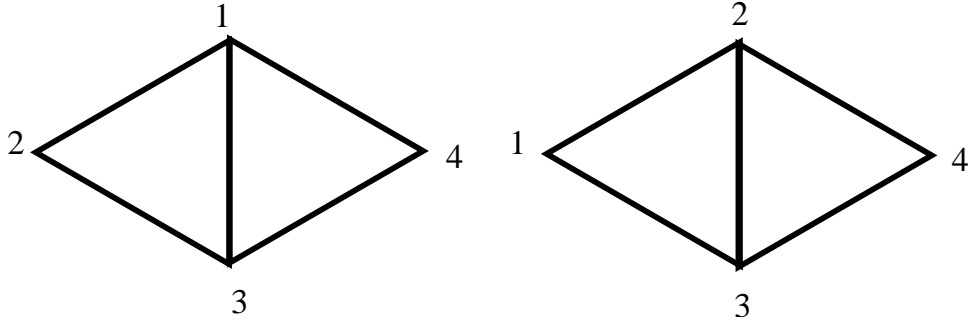}
\caption{Renaming the vertices in $\Delta$ can change the Markov complexity}
\end{center}
\end{figure}

\begin{ex}
Let $\Delta = \{\{1,2\},\{1,3\}, \{1,4\}, \{2,3\}, \{3,4\}\}$.  Then for $d_2 = d_3 = d_4 = 2$ the Markov complexity  is  $m(\Delta; 2,2,2) =  2$ while $g( \Delta; 2,2,2) = 4$.  On the other hand, for the same complex in a different orientation (see Figure 1) $\Delta' = \{\{1,2\}, \{1,3\},\{2,3\}, \{2,4\},\{3,4\}\}$ and $d_2 = d_3 = d_4 = 2$, the Markov complexity is $m( \Delta'; 2,2,2) = 4$ while $g( \Delta'; 2,2,2) = 16$. 
\end{ex}

The Graver complexity $g(A,B)$ gives an upper bound for the Markov
complexity $m(A,B)$, in terms of the Graver basis of $B$ times the
Graver basis of $A$.  There is an analogous lower bound for the
Markov complexity in terms of the Graver basis of $B$ times the
Markov basis of $A$.  To describe this lower bound, we introduce the notion of semiconformal decompositions.

\begin{defn}
Let $\bfu$, $\bfv$, and $\bfv'$ be nonzero vectors in $\ker_\zz(A)$.  We say that $\bfu
= \bfv + \bfv'$ is a \emph{semi-conformal decomposition} if $v_i > 0$
implies that $v_i \leq u_i$ and ${v'}_i < 0$ implies $u_i \leq {v'}_i$ for all $1 \leq i \leq n$.  
Note that if the first condition holds ($v_i >0$ implies that $v_i \leq u_i$ for all $i$)  the second condition ( ${v'}_i < 0$ implies $u_i \leq {v'}_i$ for all $i$) is satisfied automatically. The set 
$\mathcal{S}(A) \subset \ker_\zz(A)$ is the set 
of vectors in $\ker_\zz(A)$ which have no semi-conformal decomposition.
\end{defn}

A useful fact about vectors in the kernel of a matrix which have no semi-conformal decomposition is that they must belong to the Markov basis.

\begin{lemma}
Let $A \in \nn^{d \times n}$ with no zero columns.   If $\mathcal{M}$ is any Markov basis of $A$, then  $\mathcal{S}(A) \subseteq \mathcal{M}$.  In particular, $\mathcal{S}(A) \subseteq \mathcal{M}(A)$.
\end{lemma}

\begin{proof}
Suppose that $\bfu \in \mathcal{S}(A)$ has no semiconformal decomposition but  there is
some Markov basis $\mathcal{M}$ of $A$
that does not contain $\bfu$.
Write $\bfu = \bfu^+ - \bfu^-$ as the difference of two nonnegative integer
vectors with disjoint support.  Note that $\bfu^+$ and $\bfu^-$ belong to the
same fiber.  Since $\mathcal{M}$ is a Markov basis for $A$ there is a
sequence of elements from $\mathcal{M}$, $\{\bfv^1, \bfv^2, \ldots, \bfv^s \}$
with $s \geq 2$
which connects $\bfu^-$ to $\bfu^+$ where intermediate summands are always nonnegative.  In other words, we can write
$$ \bfu^+ = \bfu^- + \sum_{k = 1}^s \bfv^k,$$
and the set of indices where $\bfv^1$ is negative is
a subset of the set of indices where $\bfu^-$ is nonzero, and $u^-_i \geq 
|v^1_i|$ for this subset of indices.  But this implies that
$$\bfu = \bfv^1 + \sum_{k = 2}^s \bfv^k $$
\noindent is a semiconformal decomposition of $\bfu$.  This contradicts
our assumption that $\bfu \in \mathcal{S}(A)$ and hence $\bfu \in \mathcal{M}$.
\end{proof}

\begin{thm} \label{thm:lower}
The Markov complexity
$m(A,B)$ is bounded below by the maximum 1-norm of any element in
the Graver basis $\mathcal{G}( B \cdot \mathcal{S}(A))$.
\end{thm}

\begin{proof}
 Let $\Gamma \in \mathcal{G}(B \cdot \mathcal{S}(A))$.  Following the proof of 
Theorem \ref{thm:graver}, $\Gamma$ translates into a vector $\bfu =  [u^1; \ldots;
  u^r]$ with each $u^i  \in \mathcal{S}(A)$.  Furthermore, we know that $\bfu$ 
lies in the Graver basis of $\Lambda(A,B,r)$.  We
wish to show that it lies in some minimal Markov basis of $\Lambda(A,B,r)$.
To do this, we show that $\bfu$ has no semiconformal decompositions.
Suppose, to the contrary that there was some semiconformal
decomposition of $\bfu$.  Since $\bfu$ is in the Graver basis of
$\Lambda(A,B,r)$, any semiconformal decomposition $\bfu = \bfv + \bfv'$ 
induces a semiconformal decomposition of (at least) one of the
vectors $u^i$.  However, this is a contradiction, because $\mathcal{S}(A)$ consists of vectors with 
no semiconformal
decompositions. Thus $r$, which is the 1-norm of $\Gamma$ is a
lower bound for the Markov complexity $m(A,B)$.
\end{proof}

\begin{ex}
Applying Theorem \ref{thm:lower} for the simplicial complex $\Delta = \{\{1,2\}, \{1,3\},\{2,3\}\}$ we have 
$m(\Delta;3,3) \geq 5$, $m(\Delta; 3,4) \geq 8$, and $m(\Delta; 3,5) \geq 12$.  It is interesting to note that in the first two of these cases our lower bound for the Markov complexity is equal to the value reported in \cite{AT, AT2}.  The Markov complexity $M(\Delta; 3,5)$ remains undetermined.
\end{ex}


\section{Generalizations and Extensions}

Given Theorem \ref{thm:main}, it is natural to ask to what extent this 
result is the best possible.  In particular, do there exist any bounds on 
the complexity of Markov basis elements if we fix $\Delta$ and allow $d_1$ and $d_2$ to vary?  The answer to this question is negative if we allow arbitrary $\Delta$.

\begin{ex}
Let $\Delta = \{\{1,2\},\{1,3\},\{2,3\}\}$ and fix $d_3 \geq 2$.  There are no bounds $m_1$ and $m_2$ such that every Markov basis element has format strictly contained in a $m_1 \times m_2 \times d_3$ table.  A well-known example of such a move of large format is shown in \cite{DS}.  Denote by $e_{ijk}$ the $3$-way table with a $1$ in the $(i,j,k)$ position and zeroes elsewhere.  For each $m > 1$  the vector
\begin{eqnarray*}
\bfu  & =  & e_{111} + e_{221} + \cdots  + e_{mm1}  \\
  &      &   + e_{122} + e_{232} + \cdots + e_{m-1m2} + e_{m12} \\
  &      &   -  e_{112}  - e_{222}  -  \cdots - e_{mm2}  \\
  &      &    - e_{121}  - e_{231} - \cdots - e_{m-1m1} - e_{m11} 
\end{eqnarray*}
belongs to the universal Markov basis $\mathcal{M}(A_{\Delta})$ that has 
format $m \times m \times 2$.
\end{ex}

On the other hand, for reducible models, we can prove more general finiteness results.  These are based on the structural theorems for building Markov bases for reducible models in \cite{DoSu,HS}.  Recall that for a simplicial complex $\Delta$, $| \Delta| = \cup_{F \in \Delta} F$ is the underlying set of $\Delta$.

\begin{defn}
A simplicial complex $\Delta$ is called  \emph{reducible} with decomposition $(\Delta_1, S, \Delta_2)$, if $S \in \Delta_1$, $S \in \Delta_2$, $|\Delta_1| \cap |\Delta_2| = S$, and $\Delta = \Delta_1 \cup \Delta_2$.
\end{defn}

If $\Delta$ is reducible we use the notation $d^1$ and $d^2$ to denote the substrings of $d$ that are indexed by $|\Delta_1|$ and $|\Delta_2|$ respectively.

\begin{ex}
The simplicial complex $\Delta = \{\{1,2\}, \{1,3\},\{2,3\}, \{2,4\}, \{3,4\} \}$ pictured in Figure 1 is reducible with $S = \{2,3\}$, $\Delta_1 = \{\{1,2\},\{1,3\},\{2,3\}\}$, and $\Delta_2 = \{\{2,3\},\{2,4\},\{3,4\}\}$.  The vectors $d^1$ and $d^2$ are $(d_1,d_2,d_3)$ and $(d_2,d_3,d_4)$, respectively.
\end{ex}

One of the main results from \cite{DoSu, HS} is the constructive version of the following.

\begin{thm}\label{thm:gens}
Let $\Delta$ be a reducible simplicial complex and let $d$ be given.  Let $l_1$ and $l_2$ be the maximum $1$-norm of any element of 
$\mathcal{M}_{\Delta_1, d^1}$ and $\mathcal{M}_{\Delta_2, d^2}$, respectively.  Then the maximum $1$-norm of any element of $\mathcal{M}_{\Delta,d}$ is  $\max\{4,l_1,l_2\}$.
\end{thm}

This allows us to deduce that reducible models have Markov bases of finite complexity as many levels vary.

\begin{cor}
Let $\Delta$ be a reducible simplicial complex with induced subcomplexes $\Delta_1$ and $\Delta_2$ and suppose that $1 \in | \Delta_1 | \setminus | \Delta_2 |$ and $2 \in |\Delta_2 | \setminus | \Delta_1 | $.  Let $d_3, \ldots, d_n$ be given.  Then there exists constants $(m_1,m_2) = m(\Delta; d_3, \ldots, d_n)$ such that every element in the universal Markov basis $\mathcal{M}_{\Delta,d}$ has format smaller than $m_1 \times m_2 \times d_3 \times \cdots \times d_n$.
\end{cor}

\begin{proof}
Restricting to $\Delta_1$ and $\Delta_2$ and allowing $d_1$ and $d_2$ to 
vary respectively, we know by Theorem \ref{thm:main} there is a bound on 
the format of Markov basis elements that appear in 
$\mathcal{M}_{\Delta_1,d^1}$ and $\mathcal{M}_{\Delta_2,d^2}$.  However, a 
bound on the format also implies that these vectors have bounded $1$-norm 
(one such bound is the Markov complexity times the largest $1$-norm of any 
element in $\mathcal{G}(A_{{\rm link}(\Delta)})$).   Applying Theorem 
\ref{thm:gens}, we deduce that every element of $\mathcal{M}_{\Delta, d}$ 
has $1$-norm bounded by some fixed constant.  But bounded $1$-norm implies 
bounded format and completes the proof. \end{proof}

Besides the condition that $\Delta$ is reducible, the crucial requirement to prove the preceding corollary was that $1$ and $2$ were not adjacent to each other in $\Delta$.  We conjecture that this property is enough to guarantee bounded Markov complexity in general.

\begin{conj}\label{conj:indep}
Suppose that $\{1,2,\ldots,j\}$ is an independent subset of the underlying 
graph of 
$\Delta$.  Then for fixed $d_{j+1}, \ldots, d_n$, there exists numbers $(m_1, \ldots, m_j) = m(\Delta; d_{j+1}, \ldots, d_n)$ such that every element in the universal Markov basis $\mathcal{M}_{\Delta,d}$ has format smaller than $m_1 \times \cdots \times m_j \times d_{j+1} \times \cdots \times d_n$.
\end{conj}

Unfortunately, we do not know if  Conjecture \ref{conj:indep} is true even in the simplest nonreducible case, namely the four-cycle $\Delta = \{\{1,3\},\{2,3\}, \{2,4\}, \{1,4\}\}$ with $j = 2$.  If a proof exists, it must depend on techniques different from those developed here, because if $j > 1$ and $\Delta$ satisfies the hypotheses of Conjecture \ref{conj:indep} there is no bound on the formats of the Graver basis elements for $A_{\Delta}$.


\section{Computations}

The following table displays computational results of the Markov
complexity and Graver complexity of all binary hierarchical models
where one of the dimensions of the tables is allowed to vary.  In
the notation of Theorem \ref{thm:main} this is the Markov
complexity $m(\Delta;2,2,2)$ and the Graver complexity
$g(\Delta;2,2,2)$.  Note that all the entries which are marked with a
star are Markov and Graver complexities which were not known, or could not have been determined, without the use of Theorem \ref{thm:main}.  All of the computations described in this section were performed 
using the toric Gr\"obner basis program {\verb"4ti2"} \cite{He}.  The second and fifth column $m$ correspond to the computed Markov complexity and the third and sixth column $g$ is the Graver complexity.  We use the bracket notation from multivariate statistics \cite{F} for denoting simplicial complexes. Thus $[12][23][34]$ represents the simplicial complex $\{\{1,2\},\{2,3\},\{3,4\}\}$.

$$
\begin{array}{|c|c|c|}
\hline \mbox{Model} & \, \, m \, \, & \, \, g \, \, \\
\hline [123][124][134][234] & 2  & 2 \\
\hline [123][124][134]  & 2  &  2 \\
\hline [123][124][234]  & 2  &  2 \\
\hline [123][124][34]  &  2  &  2 \\
\hline [123][234][14]  &  4  &  4 \\
\hline [123][14][24][34] &  4*  &  4*  \\
\hline [234][12][13][14]  &   & 12*  \\
\hline [12][13][14][23][24][34] &   & 10*  \\
\hline [123][124]  & 2  &  2  \\
\hline [123][234]  &  2 &  2 \\
\hline [123][24][34] & 4 &  4*  \\
\hline [234][12][13] & 2 &  10* \\
\hline [123][14][24] &  2 &  2*  \\
\hline [12][13][23][24][34] & 4 & 16*  \\
\hline [12][13][14][23][24] & 2 &  4* \\
\hline [123][34]   & 2 & 2 \\
\hline [123][14]  & 2  & 2 \\
\hline [234][12]  & 2  & 4  \\
\hline  \end{array}
\begin{array}{|c|c|c|}
\hline \mbox{Model} & \, \,  m \, \, & \, \, g \, \,   \\
\hline [12][13][23][34] & 2 & 10*  \\
\hline [12][13][23][14] & 2 & 2*  \\
\hline [12][23][24][34] & 4 & 8*  \\
\hline [12][14][23][34] & 4* & 5*  \\
\hline [123][4]   & 2  & 2  \\
\hline [234][1]   & 2  & 8  \\
\hline [12][13][23][4] & 2 & 8* \\
\hline [23][24][34][1] & 4 & 16* \\
\hline [12][23][34]  & 2 & 6* \\
\hline [12][14][23]  & 2 & 4* \\
\hline [12][23][4]  & 2  & 6* \\
\hline [12][13][4]  & 2  & 3* \\
\hline [23][34][1]  & 2  & 14* \\
\hline [12][34]    & 2  & 4 \\
\hline [12][3][4]  & 2  & 4* \\
\hline [34][1][2]  & 2  & 12* \\
\hline [1][2][3][4] & 2 & 10* \\
\hline  &   &  \\
\hline
\end{array}
$$

Finally, we used the Theorem \ref{thm:lower} to compute some lower
bounds for the Markov complexity.  

$$m([12][13][23],3,5) \geq 12$$
$$m([12][13][23],4,4) \geq 16$$
$$m([123][124][134][234],3,3,3) \geq 19$$

These lower bounds are benchmark
values for extending the types of results pursued in \cite{AT} and
\cite{AT2} in which the values $m([12][13][23],3,3) = 5$ and
$m([12][13][23],3,4) = 8$ were explicitly computed.

\end{document}